\documentclass{amsart}
\usepackage{amssymb,amsmath}
\headsep=1cm \numberwithin{equation}{section}
\theoremstyle{definition}

\theoremstyle{remark}

\newcommand{\Tor}{\operatorname{Tor}}
\newcommand{\Hom}{\operatorname{Hom}}
\newcommand{\Ext}{\operatorname{Ext}}

\newcommand{\Spec}{\operatorname{Spec}}
\newcommand{\Min}{\operatorname{Min}}
\newcommand{\Max}{\operatorname{Max}}
\newcommand{\Ass}{\operatorname{Ass}}

\newcommand{\Att}{\operatorname{Att}}

\newcommand{\Supp}{\operatorname{Supp}}

\newcommand{\lrt}{\longrightarrow}

\newcommand{\Z}{\mathbb{Z}}

\newcommand{\fm}{\mathfrak{m}}
\newcommand{\fp}{\mathfrak{p}}
\newcommand{\fq}{\mathfrak{q}}

\newcommand{\fa}{\mathfrak{a}}

\newcommand{\LH}{\operatorname{H}}
\newcommand{\cd}{\operatorname{cd}}

\newcommand{\height}{\operatorname{height}}

\begin{document}

\author[M. T. Dibaei and A. Nazari]{Mohammad T. Dibaei and
 Alireza Nazari}
\title[Graded Local Cohomology]
 { Graded Local Cohomology:  attached and \\associated primes,
 asymptotic behaviors}


 \keywords{graded local cohomology, attached primes, associated primes,
  asymptotic behaviors, tameness.\\
  \hspace*{0.3cm}
 The research of the first author was in part supported from IPM (No. 00000000).\\
  Teacher Training University, Tehran, Iran and
  Institute for Studies in Theoretical Physics and Mathematics, Tehran, Iran.\\
dibaeimt@ipm.ir\\
nazaria@ipm.ir}
\begin{abstract} Assume that $R = \bigoplus_{i\in \mathbb{N}_{0}}
R_{i}$ is a homogeneous graded Noetherian ring, and that $M$ is a
$ \mathbb{Z}$--graded $R$--module, where $ \mathbb{N}_{0}$ (resp.
$ \mathbb{Z}$) denote the set all non--negative integers (resp.
integers). The set of all homogeneous attached prime ideals of the
top non--vanishing local cohomology module of a finitely generated
module $M$, $\LH _{R_{+}}^{c}(M)$, with respect to the irrelevant
ideal $R_{+}: =\bigoplus_{i\geq 1} R_{i}$ and the set of
associated primes of $\LH _{R_{+}}^{i}(M)$ is studied. The
asymptotic behavior of $\Hom_{R}(R/R_{+}, \LH _{R_{+}}^{s}(M))$
for $s \geq f(M)$ is discussed, where $f(M)$ is the finiteness
dimension of $M$. It is shown that $\LH _{R_{+}}^{h}(M)$ is tame
if $\LH _{R_{+}}^{i}(M)$ is Artinian for all $i > h$.

\end{abstract}

\maketitle

\section{Introduction}
Throughout $R = \bigoplus_{i\in \mathbb{N}_{0}} R_{i}$ is a
homogeneous positively graded Noetherian ring, so that $R =
 R_{0}[x_{1}, ..., x_{t}]$ for some $x_{1}, ..., x_{t}\in  R_{1}$,
$R_{+}: = \bigoplus_{i\geq 1} R_{i}$ is the irrelevant ideal of
$R$, and $M = \bigoplus_{i\in  \Z} M_{i}$ is a graded $R$--module
which is finitely generated whenever it is explicitly stated.
Denote by $\LH _{R_{+}}^{i}(M)$ the $i$th local cohomology module
of $M$ with respect to $R_{+}$.
 It is well--known that $\LH _{R_{+}}^{i}(M)$ inherits a
  natural grading and its $n$th components $\LH _{R_{+}}^{i}(M)_{n}$
  is finitely generated $R_{0}$--module for all $n$ and is zero
  for all $n\gg {0}$.

   Set $c : = \cd(M) = \sup \{i\in  \Z\mid (\LH
_{R_{+}}^{i}(M))\neq {0}\}$ to be the cohomological dimension of a
finitely generated module $M$ with respect to $R_{+}$. In section 2
we study  $^{*}\!\!\Att _{R}(\LH _{R_{+}}^{c}(M))$, the set all
homogeneous attached prime ideals of $\LH _{R_{+}}^{c}(M)$ and we
will show that if ($R_0, \fm_0$) is local, then the set of maximal
elements of it, $\Max (^{*}\!\Att _{R}(\LH _{R_{+}}^{c}(M)))$, is a
finite set which is equal to $\Att_{R}(\LH
_{R_{+}}^{c}(M/\fm_{0}M))$.

In section 3, we study $\Ass_{R}(\LH_{R_{+}}^{s}(M))$, for $s\geq
0$. We first show that if $K$ is a module (not necessarily finite)
over a Noetherian ring, then for each ideal $\fa$ of $S$ and each
non-negative integer $s$, we have
\begin{center}$\Ass_{S}(\LH_{\fa}^{s}(K)) \subseteq \underset {0\leq j
<s}{\cup}\Ass_{S}(\Ext_{S}^{s-j+1}(S/\fa,
\LH_{\fa}^{j}(K))/L_{j})\cup \Ass_{S}(\Ext_{S}^{s}(S/\fa,
K)/L)$\end{center} for some submodules $L, L_{0}, \cdots, L_{s-1}$
of their appropriate modules. When $M$ is a graded $R$--module
($R=R_{0}[R_{1}]$ is, as usual, a homogeneous Noetherian graded
ring) and the modules $\Ext_{R}^{s}(R/R_{+}, M)$,
$\Ext_{R}^{s-j+1}(R/R_{+}, \LH_{R_{+}}^{j}(M)), 0\leq j <s$, are $
^{\ast}$weakly Laskerian (see Definition 3.3), then
$\Ass_{R}(\LH_{R_{+}}^{s}(M))$ is a finite set.

Section 4 is devoted to asymptotic behavior of the modules $\LH
_{R_{+}}^{i}(M)$ and $\Hom_{R}(R/R_{+}, \LH _{R_{+}}^{i}(M))$. We
examine the asymptotic stability of associated primes, the
asymptotic stability of supports, and tameness of $\Hom_{R}(R/R_{+},
\LH _{R_{+}}^{i}(M))$ for graded module $M$. This gives a
generalization of the result that $\LH _{R_{+}}^{f}(M)$ is tame,
where $f: = f(M) = \inf \{i| \LH _{R_{+}}^{i}(M)$  is not finitely
generated $\}$ (cf. [BRS, Theorem 3.6(a)] ). It is shown in [Br,
Theorem 4.8 (e)] that $\LH_{R_{+}}^{\dim _{R}(M)}(M)$ is tame. We
show that $\LH_{R_{+}}^{h}(M)$ is tame if $\LH_{R_{+}}^{i}(M)$ is
Artinian $R$--module for all $i > h$. This, in particular, implies
that $\LH_{R_{+}}^{\dim _{R}(M)-1}(M)$ is tame.
\section{Attached Primes}
 Let $S$ be a commutative ring and $K$ be
an $S$--module. A prime ideal $\fp$ of $S$ is called an attached
prime of $K$ whenever there exists a submodule $L$ of $K$ such that
$\fp = L:_{S}K$ (c.f. [MS]). Note that when $K$ has a secondary
representation $K = K_{1} + ... + K_{r}$,  say,
 then $\fp = \sqrt{0:K_{i}}$ for some $i\in  \{1, ..., r\}$.  In this case
the set of attached primes of $K$ is uniquely defined (c.f. [BS,
7.2]). For the graded module $M$, if $\fp$ is an attached prime of
$M$ then it is straightforward to see that $\fp^{*}$, the
homogeneous prime ideal generated by all homogeneous elements of
$\fp$,  is also an attached prime of $M$. As usual we denote by $
^{*}\!\Spec R$
 the set of all homogeneous prime ideal of $R$. For a subset $T$
of $\Spec R$, we denote $^{{*}} T = T \cap  {^{*}\!\Spec R}$. Note
that if $M$ is an Artinian graded $R$--module, then $^{*}\!\!\Att
_{R}(M) = \Att _{R}(M)$ (cf. [Ri, corollary 1.6] ).\
  For the graded $R$--module $M$, denote by $c := \cd (M)$, the
  cohomological dimension of $M$ with respect to $R_{+}$, i.e.
  $c = \sup \{i|\LH _{R_{+}}^{i}(M) \neq {0} \}$.

  When $M$ is a finitely generated $R$--module it is shown, in [BH,
  6.2.7 ] and [Br, Proposition 3.4.(a) ], that
$c = \sup \{ \dim _{R}(M/\fm _{0}M)
  | \fm _{0}$ is a maximal element of
  $\Spec R_{0} \}$.\\

In this section, we assume that $M$ is a finitely generated graded
  $R$--module and that $R_{0}$ is not necessarily local and study the
  set $^{*}\!\!\Att _{R}(\LH
_{R_{+}}^{c}(M))$. We will show that the set $\Max (^{*}\!\Att
_{R}(\LH _{R_{+}}^{c}(M)))$, of maximal elements of $^{*}\!\!\Att
_{R}(\LH _{R_{+}}^{c}(M))$ is a finite set if $R_0$ is local (see
Corollary 2.2).
\subsection{Theorem }
{\it Assume that $M$ is a finitely generated $R$--module. Then}
\begin{center}$ \underset {\fp \in  ^{*}\!\Att _{R}(\LH _{R_{+}}^{c}(M))}{\cup}\fp =
\underset {\fp \in   \Att _{R}(\LH _{R_{+}}^{c}(M/\fm _{0}M)),
\fm_{0} \in  \Max (R_{0})}{\cup} \fp.$\end{center}

{\it Proof}. We prove it in three steps.

 {\it Step 1}. We show that $\underset{{\fp \in   \Att _{R}(\LH _{R_{+}}^{c}(M/\fm _{0}M)), \fm_{0}
\in  \Max (R_{0})}}{\cup}
  \fp\subseteq \underset{{\fp \in ^{*}\!\Att
_{R}(\LH _{R_{+}}^{c}(M))}}{\cup} \fp.$
 Note that, by
definition of $c$, for any $\fm_{0}\in \Max(R_{0})$, $\LH
_{R_{+}}^{c}(M/\fm _{0}M)$ is either zero or an Artinian graded
$R$--module and thus, by [Ri, Corollary 1.6], each element of
$\Att _{R}(\LH _{R_{+}}^{c}(M/\fm _{0}M))$, if there exists any,
is a homogeneous ideal. Choose $\fm_{0}\in \Max(R_{0})$ with $c =
\dim _{R}(M/\fm _{0}M)$. By the exact sequence $0\lrt \fm_{0}M\lrt
M\lrt M/\fm _{0}M\lrt 0$ we get an epimorphism $\LH
_{R_{+}}^{c}(M)\twoheadrightarrow \LH _{R_{+}}^{c}(M/\fm _{0}M)$
and so $\Att _{R}(\LH _{R_{+}}^{c}(M/\fm _{0}M)) \subseteq
{^{*}\!\!\Att} _{R}(\LH _{R_{+}}^{c}(M)).$ Thus our claim is
clear.

 {\it Step 2}. We prove that $\underset{\fp \in  ^{*}\!\!\Att _{R}(\LH
_{R_{+}}^{c}(M))}{\cup} \fp \subseteq \underset{\fp \in
^{*}\!\Supp_{R}(M), \cd(R/\fp) = c }{\cup} \fp.$ Set $X = \{\fp \in
\Ass (M) |\cd(R/\fp)\\
 = c\}$. Thus there exists a graded submodule
$N$ of $M$ such that $\Ass_{R}(M/N) = X$ and $\Ass_{R}(N) =
\Ass_{R}(M)\setminus X$ ( one may choose $N$ to be a maximal element
of the set $\{T | T$ is a graded submodule of $M$ such that
$\Ass_RT\subseteq \Ass_R(M)\setminus X\}$ (see [Bo, page 263,
Proposition 4] for non--graded case) ). Consider the exact sequence
$\LH _{R_{+}}^{c}(N)\lrt \LH _{R_{+}}^{c}(M) \lrt \LH
_{R_{+}}^{c}(M/N)\lrt 0$ of graded modules. Note that $\cd(N)\leq c$
(c.f. [Br, Corollary 3.5]) and $\Ass_{R}(N)\subseteq
\Ass_{R}(M)$.Therefore it follows that $\LH _{R_{+}}^{c}(N)=0$ and
$\LH _{R_{+}}^{c}(M)\cong \LH _{R_{+}}^{c}(M/N)$. Now, choose $r$ to
be a homogeneous element of $R$ such that $r\not\in \underset{\fp
\in ^{*}\!\Supp_{R}(M), \cd(R/\fp) = c }{\cup} \fp$. It follows from
definition of $N$ that $r$ is not a zero divisor on $M/N$. Consider
the exact sequence $0\lrt M/N[-l] \overset{r}\lrt M/N\lrt
\frac{M/N}{r(M/N)} \lrt 0,   (l := \deg r),$ from which we have the
exact sequence $ \LH _{R_{+}}^{c}(M/N)[-l]\overset{r}\lrt \LH
_{R_{+}}^{c}(M/N)\lrt \LH _{R_{+}}^{c}(\frac{M/N}{r(M/N)}) \lrt 0.$
If $\LH _{R_{+}}^{c}(\frac{M/N}{r(M/N)})\neq0$, then $\cd(\frac
{M/N}{r(M/N)}) = c$ and so $\cd (R/\fq) = c$ for some $q\in
\Ass_{R}(\frac{M/N}{r(M/N)})$. As $q\in {^{*}\!\Supp_{R}(M)}$ and
$r\in  \fq$, this gives a contradiction. Therefore we have  $\LH
_{R_{+}}^{c}(M/N) = r \LH _{R_{+}}^{c}(M/N)$ which implies that $r
\not\in  {\cup}_{\fp \in ^{*}\!\!\Att _{R}(\LH _{R_{+}}^{c}(M/N))}
\fp$. This completes step 2 because $\LH _{R_{+}}^{c}(M)\cong \LH
_{R_{+}}^{c}(M/N)$.

 {\it Step 3}.
Finally we show $\underset{\cd(R/\fp) = c, \fp \in
^{*}\!\Supp_{R}(M)}{\cup} \fp \subseteq  \underset {\fm_{0} \in \Max
(R_{0}), \fp \in \Att _{R}(\LH _{R_{+}}^{c}(M/\fm _{0}M))}{\cup}{
\fp}$. Choose $\fp \in {^{*}\!\Supp_{R}(M)}$ with $\cd(R/\fp) = c$ .
Hence $\dim _{R}(\frac{R/\fp}{\fm_{0}(R/\fp)}) = c$ for some
$\fm_{0}\in \Max(R_{0})$ and so there exists $\fq\in \Supp_{R}(M/\fm
_{0}M)$ such that $\fm_{0}R + \fp \subseteq \fq$ and $\dim
_{R}(R/\fq) = c$. This implies that $\fq\in \Att _{R}(\LH
_{R_{+}}^{c}(M/\fm _{0}M))$ (see [MS], or [DY1, Theorem A] for
non--local case). This is step 3 and so the proof is complete.
$\hfill\Box$\\

The following corollary is in contrast to the result [Br,
Corollary 3.9] which states that $\Min(\Ass _{R}(\LH
_{R_{+}}^{c}(M)))$ is a finite set.
\subsection{Corollary }
{\it Assume that $(R_{0}, \fm_{0})$ is a local ring and that $M$ is
a finitely generated graded $R$--module. Then}
$$\Max(^{*}\!\Att _{R}(\LH _{R_{+}}^{c}(M)))
 = \Att _{R}(\LH_{R_{+}}^{c}(M/\fm _{0}M)).$$ {\it In particular,}
$\Max(^{*}\!\Att _{R}(\LH _{R_{+}}^{c}(M)))$ {\it is a finite set
and it depends only on} $\Supp_{R}(M).$

 {\it Proof}. Note that the set $\Att _{R}(\LH_{R_{+}}^{c}(M/\fm _{0}M))$ is a subset of
 the associated height of $M/\fm _{0}M$ so that it is a finite set
 and there are no containment relations among its elements
  (c.f. [DY1, Theorem A]).
 Hence, by Theorem 2.1, and the right exactness of $\LH_{R_{+}}^{c}(-)$,
 we get the result.  $\hfill\Box$\\

\indent Now we present an example of a finitely generated graded
module $M$ such that $^{*}\!\Att _{R}(\LH _{R_{+}}^{c}(M))$ is not a
finite set. Assume that $R=R_0[X]$ is the polynomial ring over a
Noetherian local ring $R_0$ with $\dim R_0>1$. As $X$ is a non--zero
divisor on $R$, it is trivial that $\LH_{(X)}^{1}(R)\cong R_{X}/R$.
It is easy to see that $\fp_0 [X]\in ^{*}\!\Att _{R}(\LH
_{(X)}^{1}(R))$ for any $\fp_0 \in \Spec R_0$ so that $^{*}\!\Att
_{R}(\LH _{(X)}^{1}(R))$ is an infinite set.

\section{Associated Primes, non--graded and graded cases}
For a graded $R$--module $M$, we study the set of associated primes
of $\LH_{R_{+}}^{s}(M)$. In this section we first bring some results
for non--graded case and then we state them for our main purpose. We
assume that $S$ is a Noetherian ring, $K$ is an $S$--module and
$\fa$ is an ideal of $S$. The aim of the following result is to show
that, for any non--negative integer $t$, the set of associated
primes of $\LH_{\fa}^{t}(K)$ depends on the set of the associated
primes of some quotients of the modules $\Ext_{S}^{t}(S/\fa, K)$ and
$\Ext_{S}^{t-j+1}(S/\fa, \LH_{\fa}^{j}(K)), 0\leq j <t$ (compare
with [DM, Theorem 2.5]).
\subsection{Theorem}
{\it Assume that $S$ is a Noetherian ring, $\fa$ is an ideal of $S$
and that $K$ is an $S$--module. Then, for each non--negative integer
$t$, there exist a submodule $L$ of $\Ext_{S}^{t}(S/\fa, K)$ and
submodules $L_{j}$ of $\Ext_{S}^{t-j+1}(S/\fa, \LH_{\fa}^{j}(K)),
0\leq j <t$, such that}
\begin{center}$\Ass_{S}(\LH_{\fa}^{t}(K)) \subseteq \cup_{0\leq j
<t}\Ass_{S}(\Ext_{S}^{t-j+1}(S/\fa, \LH_{\fa}^{j}(K))/L_{j})\cup
\Ass_{S}(\Ext_{S}^{t}(S/\fa, K)/L).$\end{center}

{\it Proof}. The result is clear for $t=0$. Assume that $t
>0$ and that $t-1$ is settled. Note that $\LH_{\fa}^{j}(K)\cong
\LH_{\fa}^{j}(K/\Gamma_{\fa}(K))$ for all $j
> 0$. Assume that $E$ is an injective hull of $K/\Gamma_{\fa}(K)$
and set $N:=E/(K/\Gamma_{\fa}(K))$. Hence we have
$\LH_{\fa}^{i}(N)\cong \LH_{\fa}^{i+1}(K)$ and $\Ext_{S}^{i}(S/\fa,
N)\cong \Ext_{S}^{i+1}(S/\fa, K/\Gamma_{\fa}(K))$ for all $i\geq 0$.
By our induction hypothesis, there exist submodules $T$ and
$L_{j+1}$ of the
appropriate modules  such that \\[-0.6 cm]\begin{small}
\begin{equation} \Ass_{S}(\LH_{\fa}^{t-1}(N)) \subseteq \underset {0\leq j
<t-1} {\cup}\Ass_{S}(\frac{\Ext_{S}^{t-1-j+1}(S/\fa,
\LH_{\fa}^{j}(N))}{L_{j+1}})\cup
\Ass_{S}(\frac{\Ext_{S}^{t-1}(S/\fa,
N)}{T}).\tag{1}\end{equation}\end{small}
 Note that from the exact sequence
$\Ext_{S}^{t}(S/\fa, K)\overset{f}{\longrightarrow}
\Ext_{S}^{t}(S/\fa, K/\Gamma_{\fa}(K))\overset{g}{\longrightarrow}
\Ext_{S}^{t+1}(S/\fa, \Gamma_{\fa}(K)) $ we have the induced exact
sequence $$\Ext_{S}^{t}(S/\fa,
K)/f^{-1}(T)\overset{\overline{f}}{\longrightarrow}
\Ext_{S}^{t}(S/\fa,
K/\Gamma_{\fa}(K))/T\overset{\overline{g}}{\longrightarrow}
\Ext_{S}^{t+1}(S/\fa, \Gamma_{\fa}(K))/g(T).$$ Therefore, there are
submodules $L$ and $L_0$ of the appropriate modules such that
\begin{small}\begin{equation} \Ass_{S}(\frac{\Ext_{S}^{t}(S/\fa,
K/\Gamma_{\fa}(K))}{T})\subseteq \Ass_{S}(\frac{\Ext_{S}^{t}(S/\fa,
K)}{L})\cup \Ass_{S}(\frac{\Ext_{S}^{t+1}(S/\fa,
\Gamma_{\fa}(K))}{L_{0}}). \tag{2}\end{equation}\end{small} Now, (1)
and (2) imply the claim. $\hfill \Box$\\
\subsection{Corollary}
{\it Assume that $R = \underset{i\in
\mathbb{N}_{0}}{\bigoplus}R_{i}$ is a homogeneous Noetherian
graded ring (i.e. $R = R_{0}[R_{1}])$ and $M$ is a graded
$R$--module. Then, for each non--negative integer $s$, there exist
graded submodules $L$ of $\Ext_{R}^{s}(R/R_{+}, M)$ and $L_{j}$ of
$\Ext_{R}^{s-j+1}(R/R_{+}, \LH_{R_{+}}^{j}(M)), 0 \leq j < s,$
such that $$\Ass_{R}(\LH_{R_{+}}^{s}(M)) \subseteq \underset
{0\leq j <s}{\cup}\Ass_{R}(\frac{\Ext_{R}^{s-j+1}(R/R_{+},
\LH_{R_{+}}^{j}(M))}{L_{j}})\cup
\Ass_{R}(\frac{\Ext_{R}^{s}(R/R_{+}, M)}{L}).$$}

{\it Proof.} The proof is similar to that of Theorem 3.1. One might
take into consideration that in the proof we replace injective hull
of $M/\Gamma_{R_{+}}(M)$ by $ ^{\ast}$injective hull of
$M/\Gamma_{R_{+}}(M)$ and note that all modules are graded and all
homomorphisms are homogeneous. $\hfill \Box$\\
\subsection{Definition}(see [DM, Definition 2.1])
A graded $R$--module $M$ is called $ ^{\ast}${\it weakly
Laskerian} if for each graded submodule $N$ of $M$,
$\Ass_{R}(M/N)$ is a finite set.
\subsection{Corollary}
{\it Assume that $R$ is a homogeneous Noetherian graded ring, $M$ is
a graded $R$--module and that $s$ is a non--negative integer. If all
the modules $\Ext_{R}^{s}(R/R_{+}, M)$ and  $$
\Ext_{R}^{s-j+1}(R/R_{+}, \LH_{R_{+}}^{j}(M)), 0 \leq j < s,$$
 are $^{\ast}$weakly Laskerian, then $\Ass_{R}(\LH_{R_{+}}^{s}(M))$ is a
finite set.}

{\it Proof.} It follows from Corollary 3.2. $\hfill \Box$\\
\section{Asymptotic behaviors}
Assume that $M$ is a graded $R$--module. The module $M$ is said to
have the property of {\it asymptotic stability of associated
primes} (resp. {\it asymptotic stability of supports}) if there
exists an integer $n_{0}$ such that $\Ass_{R_{0}}(M_{n}) =
\Ass_{R_{0}}(M_{n_{0}})$ for all $n\leq n_{0}$ (resp.
$\Supp_{R_{0}}(M_{n}) = \Supp_{R_{0}}(M_{n_{0}})$ for all $n\leq
n_{0}$). The module $M$ is called {\it tame} if $M_{i}=0$ for all
$i\ll 0$ or else $M_{i}\neq 0$ for all $i\ll 0$. Here, we study
the above asymptotic behaviors of $\Hom_{R}(R/R_{+}, \LH
_{R_{+}}^{i}(M))$. In this connection, there are three open
problems.
\subsection{Problems}
 (cf. [Br, Problems 6.1, Problem 7.1, and
Problem 4.3]) Let $i\in \mathbb{N}_{0}$ and let $M$ be a finitely
generated graded $R$--module.\\
\indent (i) $\LH_{R_{+}} ^{i}(M)$ has the property of asymptotic
stability of associated primes.\\
\indent (ii) $\LH_{R_{+}} ^{i}(M)$ has the property of asymptotic
stability of supports.\\
\indent (iii) $\LH_{R_{+}} ^{i}(M)$ is tame.\\

Note that (i) implies (ii) and (ii) implies (iii). In this section
we investigate the above questions for the module
$\Hom_{R}(R/R_{+}, \LH _{R_{+}}^{i}(M))$. We note that for a
graded $R$--module $M$, $\Ext_{R}^{i}(R/R_{+}, M)$ has a natural
grading and $\Ext_{R}^{i}(R/R_{+}, M)\cong
{^{^{*}}\!\!\Ext_{R}^{i}(R/R_{+}, M)}$ for all $i\geq0$ (see [BS,
Proposition 12.2.7]). We first note the following easy lemma.
\subsection{Lemma}
{\it Assume that $R$ is a homogeneous Noetherian graded ring and
that $M$ is
a graded $R$--module. Then the following statements hold.\\
(i). If $\Hom_{R}(R/R_{+}, M)$ is tame, then $M$ is tame.\\
(ii). Assume that ($R_{0}, \fm_{0}$) is local and that $\fm^{\ast}:=
\fm_{0} + R_{+}$ is  the $ ^{\ast}$maximal ideal of $R$. Assume
that, for the graded $R$--module $M$, each component $M_{i}$ is a
finite $R_{0}$--module. If $M/\fm^{\ast}M$ is Artinian, then $M$ is
tame.}

{\it Proof.} (i). Assume that there is $n_{0}\in  \mathbb Z$ such
that $\Hom_{R}(R/R_{+}, M)_{n}$ is either zero for all $n$ with
$n\leq n_{0}$ or is non--zero for all $n$ with $n\leq n_{0}.$ Now
assume that $M_{k} = 0$ for some $k\leq n_{0}.$ We show that $M_{l}
= 0$ for all $l\leq k.$ It is enough to show that $M_{k-1} = 0.$ As
$R_{1}M_{k-1}\subseteq M_{k} = 0,$ we have $ {M_{k-1}} \subseteq
0\underset{M_{k-1}}{:}R_{+} \cong\Hom_{R}(R/R_{+}, M) _{k-1} =0. $\\
\indent (ii). As $M/\fm^{\ast}M$ is Artinian  it is also Noetherian
and so $(M/\fm^{\ast}M)_{k} = 0$ for all $k\ll 0$. By Nakayama Lemma
we get $M_{k}=R_{1}M_{k-1}$ for all $k\ll 0$. This implies that $M$
is tame. $\hfill \Box$\\
\subsection{Definition}
A graded module $M$ over a homogeneous graded ring $R$ is called
{\it asymptotically zero} if $M_{n}=0$ for all $n\ll 0$.

All finitely generated graded $R$--modules are asymptotically
zero. Now, we are ready to present our main results of this
section. We put these results in the following theorem.
\subsection{Theorem}
{\it Assume that $M$ is a graded $R$--module and that $s$ is a fixed
non--negative integer such that the modules
$$\Ext_{R}^{s-j}(R/R_{+}, \LH_{R_{+}}^{j}(M)),
\Ext_{R}^{s-j+1}(R/R_{+}, \LH_{R_{+}}^{j}(M)), j=0,1,\cdots, s-1 $$
are asymptotically zero (e.g. they might be finitely generated).
Then the following statements hold.

(i) The module $\Ext_{R}^{s}(R/R_{+}, M)$ has the property of
asymptotic stability of associated primes if and only if
$\Hom_{R}(R/R_{+}, \LH_{R_{+}}^{s}(M))$ has the property of
asymptotic stability of associated primes.

(ii) The module $\Ext_{R}^{s}(R/R_{+}, M)$ has the property of
asymptotic stability of supports if and only if $\Hom_{R}(R/R_{+},
\LH_{R_{+}}^{s}(M))$ has the property of asymptotic stability of
supports.

(iii)  The module $\Ext_{R}^{s}(R/R_{+}, M)$ is tame if and only
if $\Hom_{R}(R/R_{+}, \LH_{R_{+}}^{s}(M))$ is tame.}

{\it Proof.} The proofs of (i), (ii), and (iii) are essentially
similar therefore we give a proof for (iii) only. The proof is
inspired by that of [DY2, Theorem 6.3.9]. The case $s$ = 0 is
trivial because we have $\Hom_{R}(R/R_{+}, \Gamma_{R_{+}}(M))=
\Hom_{R}(R/R_{+}, M)$ . Assume that $s>0$ and that the case $s-1$ is
settled. Denote the $ ^{^{*}}\!\!$ injective hull of
$M/\Gamma_{R_{+}}(M)$ by $^{^{*}}\!\!E$ and denote
$N:={^{^{*}}\!\!E}/(M/\Gamma_{R_{+}}(M)).$ Therefore we have the
exact sequence $$0\longrightarrow M/\Gamma_{R_{+}}(M)\longrightarrow
{^{^{*}}\!\!E} \longrightarrow N \longrightarrow 0, $$ which implies
the isomorphisms $\LH_{R_{+}}^{j}(N)\cong \LH_{R_{+}}^{j+1}(M)$ and
\begin{equation} \Ext_{R}^{j}(R/R_{+}, N)\cong \Ext_{R}^{j+1}(R/R_{+},
M/\Gamma_{R_{+}}(M)) \tag{\dag} \end{equation} for all $j\geq 0$
(note that $\Hom_{R}(R/R_{+}, {^{^{*}}\!\!E}) = 0 = \Gamma_{R_{+}}(
{^{^{*}}\!E}))$. Thus, for all $j\geq 0$, we have the isomorphisms
\begin{center}$\Ext_{R}^{s-1-j}(R/R_{+}, \LH_{R_{+}}^{j}(N))\cong
\Ext_{R}^{s-(j+1)}(R/R_{+}, \LH_{R_{+}}^{j+1}(M))$\end{center} and
\begin{center}$\Ext_{R}^{s-1-j+1}(R/R_{+}, \LH_{R_{+}}^{j}(N))\cong
\Ext_{R}^{s-(j+1)+1}(R/R_{+}, \LH_{R_{+}}^{j+1}(M)).$\end{center}
Therefore, for all $j=0, \cdots, s-2$, the modules
\begin{center}$\Ext_{R}^{s-1-j}(R/R_{+}, \LH_{R_{+}}^{j}(N)) ,
\Ext_{R}^{s-1-j+1}(R/R_{+}, \LH_{R_{+}}^{j}(N))$\end{center} are
asymptotically zero. We now prepare the requirement for the
induction step $s-1$ for $N$. By the exact sequence
$0\longrightarrow \Gamma_{R_{+}}(M)\longrightarrow M
\longrightarrow M/\Gamma_{R_{+}}(M)\longrightarrow 0$ we have the
exact sequence of graded modules with homogeneous homomorphisms
\begin{center}
$\Ext_{R}^{s}(\frac{R}{R_{+}}, \Gamma_{R_{+}}(M))\rightarrow
\Ext_{R}^{s}(\frac{R}{R_{+}}, M)\rightarrow
\Ext_{R}^{s}(\frac{R}{R_{+}},
\frac{M}{\Gamma_{R_{+}}(M)})\rightarrow
\Ext_{R}^{s+1}(\frac{R}{R_{+}}, \Gamma_{R_{+}}(M)).$
\end{center}
 By our hypothesis for $s$ and $j=0$, there exists $t \in  \Z$ such that
the modules $\Ext_{R}^{s}(R/R_{+}, \Gamma_{R_{+}}(M))_{n}$ and
$\Ext_{R}^{s+1}(R/R_{+}, \Gamma_{R_{+}}(M))_{n}$ are zero for all
$n\leq t$. Therefore, for each $n\leq t$, we have an
$R_{0}$--isomorphism \\[-0.8 cm]
\begin{equation}\Ext_{R}^{s}(R/R_{+}, M)_{n} \cong \Ext_{R}^{s}(R/R_{+},
M/{\Gamma_{R_{+}}(M)})_{n}.\tag{\dag\dag}\end{equation}\\[-0.8 cm]
 Now, it follows, from ($\dag$) and ($\dag\dag$), that  $\Ext_{R}^{s}
 (R/R_{+}, M)$ is tame if and only if
 $\Ext_{R}^{s-1}(R/R_{+}, N)$ is tame, which is also equivalent to say
 that $\Hom_{R}(R/R_{+},
\LH_{R_{+}}^{s-1}(N))$ is tame, by our induction hypothesis. This
statement is also equivalent to say $\Hom_{R}(R/R_{+},
\LH_{R_{+}}^{s}(M))$ is tame.   $\hfill\Box$\\
\subsection{Corollary}
(see [Br, Theorem 4.8(b)] and [DY3, Theorem 2.1]) {\it Assume that
$M$ is a finitely generated graded $R$--module and that $s$ is a
fixed non--negative integer such that for each $i<s$,
$\LH_{R_{+}}^{i}(M)$ is $R_{+}$--cofinite, i.e.
$\Ext_{R}^{j}(R/R_{+}, \LH_{R_{+}}^{i}(M))$ is finitely generated
for all $j$. Then $\Hom_{R}( R /R_{+}, \LH_{R_{+}}^{s}(M))$ is
asymptotically zero. In particular, $\Hom_{R}( R /R_{+},
\LH_{R_{+}}^{s}(M))$ is finitely generated and so
$\LH_{R_{+}}^{s}(M)$ is tame.} $\hfill\Box$\\

Let $M$ be a finitely generated $R$--module and $n: =\dim _{R}(M)$.
In [Br, Theorem 4.8 (e)], Brodmann showed that $\LH_{R_{+}}^{c}(M)$
is tame. This clearly implies that $\LH_{R_{+}}^{n}(M)$ is tame
(this is also a consequence of the fact that $\LH_{R_{+}}^{n}(M)$ is
an Artinian module). Brodmann also showed that $\LH_{R_{+}}^{i}(M)$
is tame for all $i\in \Psi(M)$ , where $\Psi(M) = \{ \height (\fp/(
0 :_{R} M)) |  \fp \in \Min(( 0:_{R}M )) + R_{+} \}$ (see [Br,
Theorem 4.8 (d)]). Note that when $R_{0}$ is a field and $M =
R_{0}[x]$, we have $\Psi(M) = \{\dim _{R}(M)\}$. Therefore it would
be significant to see explicitly that $\LH_{R_{+}}^{n-1}(M)$ is
tame. This result is a consequence of a more general one ( see
Corollary 4.8).

In the remainder of this section we assume that $M$ is a graded
module over the homogeneous Noetherian graded ring $R =
\underset{i\in \mathbb{N}_{0}}{\bigoplus}R_{i}$.
\subsection{Definition}(See [DY3, Definition 3.1])
Define $q(M) := \sup \{i | H_{R_{+}}^{i}(M)$ is not Artinian$\}.$
If $H_{R_{+}}^{i}(M)$ is Artinian for all $i$, we write $q(M)=
-\infty$.
\subsection{Theorem}
{\it Assume that the base ring $R_0$ is local with maximal ideal
$\fm_0$ and that $M$ is a finitely generated graded $R$--module.
Then
$$\LH_{R_{+}}^{q(M)}(M)/\fm^{\ast}\LH_{R_{+}}^{q(M)}(M)$$ is
Artinian.}

{\it Proof.}  Set $n(M):=\cd (M) - q(M)$ and we prove our claim by
using induction on $n(M)$. When $n(M)=0$ the result is known (c.f.
[Br, Theorem 2.3 (b)]) because, in this case,
$\LH_{R_{+}}^{q(M)}(M)/\fm^{\ast}\LH_{R_{+}}^{q(M)}(M)$ is a
homomorphic image of
$\LH_{R_{+}}^{q(M)}(M)/\fm_{0}\LH_{R_{+}}^{q(M)}(M)$ . Assume that
$n(M)=n
>0$ and we have proved the statement for any finitely generated
graded $R$--module $N$ with $n(N)= n-1$. Thus we have $\cd (M) > 0$.
We may assume that $q(M)> 0$. Therefore, $q(M) =
q(M/\Gamma_{R_{+}}(M))$, $\cd (M) = \cd (M/\Gamma_{R_{+}}(M))$, and
$n(M) = n(M/\Gamma_{R_{+}}(M))$. Hence we may assume that $M$ is
$R_{+}$--torsion free and there exists $x\in R_{1}$ which is
non-zero-divisor on $M$ and $\cd (M/xM) = \cd (M)-1$ (c.f. [RS,
1.3.7]). The short exact sequence $0\longrightarrow M[-1]
\overset{x}{\longrightarrow}M \longrightarrow M/xM \longrightarrow
0$ yields a long exact sequence
$$\cdots \longrightarrow
\LH_{R_{+}}^{q(M)}(M)[-1]\overset{x}{\longrightarrow}\LH_{R_{+}}^{q(M)}(M)
\longrightarrow \LH_{R_{+}}^{q(M)}(M/xM) \longrightarrow
\LH_{R_{+}}^{q(M)+1}(M)[-1] \longrightarrow \cdots $$ from which we
have the exact sequence $$ 0 \longrightarrow
\LH_{R_{+}}^{q(M)}(M)/x\LH_{R_{+}}^{q(M)}(M) \longrightarrow
\LH_{R_{+}}^{q(M)}(M/xM) \longrightarrow ( 0
\underset{\LH_{R_{+}}^{q(M)+1}(M)}{:} x )[-1] \longrightarrow 0.$$
Therefore we have the exact sequence $$\Tor_{1}^{R} (R/\fm^{\ast}, (
0 \underset{\LH_{R_{+}}^{q(M)+1}(M)}{:} x))[-1] \longrightarrow
\LH_{R_{+}}^{q(M)}(M)/\fm^{\ast}\LH_{R_{+}}^{q(M)}(M)
\longrightarrow$$ $$
\LH_{R_{+}}^{q(M)}(M/xM)/\fm^{\ast}\LH_{R_{+}}^{q(M)}(M/xM)
 \longrightarrow ( 0 \underset{\LH_{R_{+}}^{q(M)+1}(M)}{:}
x)/\fm^{\ast}( 0 \underset{\LH_{R_{+}}^{q(M)+1}(M)}{:} x)[-1].$$
Note that the first and the last term in the above exact sequence
are Artinian modules. If $\LH_{R_{+}}^{q(M)}(M/xM)$ is Artinian then
$\LH_{R_{+}}^{q(M)}(M/xM)/\fm^{\ast}\LH_{R_{+}}^{q(M)}(M/xM)$ is
Artinian and so
$$\LH_{R_{+}}^{q(M)}(M)/\fm^{\ast}\LH_{R_{+}}^{q(M)}(M)$$ is also
Artinian.  Now, we assume that $\LH_{R_{+}}^{q(M)}(M/xM)$ is not
Artinian. It follows that $q(M/xM) = q(M)$ and hence $n(M/xM) = \cd
(M/xM)-q(M/xM)=n(M)-1$. By our induction hypothesis,
$\LH_{R_{+}}^{q(M)}(M/xM)/\fm^{\ast}\LH_{R_{+}}^{q(M)}(M/xM)$ is
Artinian. By the above exact sequence, the module
$$\LH_{R_{+}}^{q(M)}(M)/\fm^{\ast}\LH_{R_{+}}^{q(M)}(M)$$ is Artinian.
$\hfill\Box$\\
\subsection{Corollary}
{\it With the assumptions as in Theorem 4.7, $\LH_{R_{+}}^{i}(M)$ is
tame for all $i\geq q(M)$. In particular, the modules
 $\LH_{R_{+}}^{\cd (M)}(M)$ and $\LH_{R_{+}}^{\dim_{R}(M)-1}(M)$ are
tame.}

{\it Proof.} Note that $\LH_{R_{+}}^{\dim_{R}(M)}(M)$ is Artinian
(c.f. [BS, Theorem 7.1.6]). The claims follow by Theorem 4.7 and
Lemma 4.2(ii). $\hfill\Box$\\

{\bf Acknowledgment} The first author would like to express his
thank to Professor J\"{u}rgen Herzog for his informative discussion
about graded modules during his visit to Mathematics Department in
University of Duisburg-Essen. The authors would like to thank the
referee for her/his comments.


\end{document}